# Modeling of Nonlinear Dynamic Systems with Volterra Polynomials: Elements of Theory and Applications

A.S. Apartsyn, S.V. Solodusha, V.A. Spiryaev

*Energy Systems Institute, Siberian Branch of the Russian Academy of Sciences, Russia*

**ABSTRACT**

The paper presents a review of the studies that were conducted at Energy Systems Institute (ESI) SB RAS in the field of mathematical modeling of nonlinear input-output dynamic systems with Volterra polynomials. The first part presents an original approach to identification of the Volterra kernels. The approach is based on setting special multi-parameter families of piecewise constant test input signals. It also includes a description of the respective software; presents illustrative calculations on the example of a reference dynamic system as well as results of computer modeling of real heat exchange processes. The second part of the review is devoted to the Volterra polynomial equations of the first kind. Studies of such equations were pioneered and have been carried out in the past decade by the laboratory of ill-posed problems at ESI SB RAS. A special focus in the paper is made on the importance of the Lambert function for the theory of these equations.

**INTRODUCTION**

In the theory of mathematical modeling of nonlinear dynamic systems the universal technique of Volterra functional series is well known (Volterra, 1982). It presents a response $y(t)$ of the input-output system to an external disturbance $x(t)$ (for simplicity consider $x(t)$ and $y(t)$ to be scalar time functions) as the integro-power series

$$y(t) = F(x(t)) \equiv \sum_{m=1}^{\infty} \int_0^t \ldots \int_0^t K_m(t, s_1, \ldots, s_m) \prod_{i=1}^m x(s_i) ds_i, \quad t \in [0, T]. \qquad (1)$$

In (1) the functions $K_m$, $m = 1,2,...$, that are called Volterra kernels, are identified on the basis of information about system responses to certain families of test input signals.

Application of the finite segment of series (1)

$$y(t) = P_N(x(t)) \equiv \sum_{m=1}^{N} \int_0^t ... \int_0^t K_m(t, s_1, ..., s_m) \prod_{i=1}^{m} x(s_i) ds_i, \quad t \in [0, T], \quad (2)$$

to modeling of nonlinear dynamic objects of various nature is based on continual analogs of the classical Weierstrass theorem about approximation of continuous function by polynomial, i.e. Frechet theorem (Frechet, 1910) and its different generalizations. For example, according to (Baesler & Daugavet, 1990) for any Volterra initial $(F(x(0)) = 0)$ mapping $y(t) = F(x(t))$, $F: \mathcal{K} \to C_{[0,T]}$, $\mathcal{K}$ – compact in $C_{[0,T]}$, and any $\varepsilon > 0$ there is $N$ such that for all $x(t) \in K$ the equality

$$\|F(x(t)) - P_N(x(t))\|_{C_{[0,T]}} < \varepsilon,$$

holds, where $P_N(x(t))$ has the form (2).

The problem of applying the Volterra series to modeling of technical systems (in particular electric power and heat supply systems) was considered in the monographs (Pupkov, Kapalin & Yushchenko, 1976; Deitch, 1979; Venikov & Sukhanov, 1982; Pupkov & Shmykova, 1982; Danilov, Matkhakov & Filippov, 1990). The international sources to be emphasized are the profound monographs (Rugh, 1981; Doule III, Pearson & Ogunnaike, 2002; Ogunfunmi, 2007) that provide an extensive bibliography.

If a system to be modeled is stationary in the sense that its dynamic characteristics during the transient process $T$ can be considered constant, the Volterra kernels $K_m$ depend only on the differences $t - s_i, i = \overline{1, m}$, and then instead of (2) the model

$$y(t) = \sum_{m=1}^{N} \int_0^t ... \int_0^t K_m(s_1, ..., s_m) \prod_{i=1}^{m} x(t - s_i) ds_i, \quad t \in [0, T], \quad (3)$$

is used in which $K_m$ are continuous and symmetric with respect to a set of variables. At $N = 1$ (3) presents the output signal through the convolution integral (Duhamel integral)

$$y(t) = \int_0^t K_1(s) x(t - s) ds, \quad t \in [0, T], \quad (4)$$

and this representation is standard for the linear theory of automatic control.

In this case for identification of $K_1(t)$ it is sufficient to know a response of the system to a random disturbing action. If, in particular, to assume

$$x(t) = \alpha I(t), \quad (5)$$

where $I(t)$ – the Heaviside function, and scalar $\alpha \neq 0$, then at $y(t) \in C^{(1)}_{[0,T]}, y(0) = 0$, the formula of inversion (4) is true

$$K_1(t) = \frac{y'(t)}{\alpha}. \tag{6}$$

## ABOUT ONE APPROACH TO IDENTIFICATION OF VOLTERRA KERNELS

The ESI SB RAS started the research in the field of the Volterra kernels identification in the 1990s. In (Apartsyn, 1991; Apartsyn, 1992), for the case $N = 2,3$, and in (Apartsyn, 1995a; Apartsyn, 1995b; Apartsyn, 1995c; Apartsyn, 1996) for arbitrary $N$ the method of identification of $K_m, m = \overline{1,N}$ was proposed. It was natural generalization of (5), (6) for a multidimensional case. For identification of $K_m$ we introduced an $(m - 1)$ – parameter family of piecewise constant functions as test signals

$$x^{\alpha_k}_{\omega_1, \ldots, \omega_{m-1}}(t) = \alpha_k \sum_{i=0}^{m-1} \gamma_i I(t - \beta_i), \tag{7}$$

where $\gamma_i = (-1)^i \cdot 2, i = \overline{1, m-2}$; $\gamma_0 = 1, \gamma_{m-1} = (-1)^{m-1}$; $\beta_i = \sum_{j=1}^{i} \omega_j$, $\beta_0 = 0$, $\omega_1, \ldots, \omega_{m-1}, t \in [0, T]$; $k = \overline{1, m}$, $\alpha_1 \neq \alpha_2 \neq \cdots \neq \alpha_m \neq 0$.

In particular for the most important practical cases $m = 2,3$

$$x^{\alpha_k}_{\omega_1}(t) = \alpha_k [I(t) - I(t - \omega_1)], k = 1,2; \tag{8}$$

$$x^{\alpha_k}_{\omega_1, \omega_2}(t) = \alpha_k [I(t) - 2I(t - \omega_1) + I(t - \omega_1 - \omega_2)], k = 1,2,3; \tag{9}$$

If to understand under $f_m(t, \omega_1, \ldots, \omega_{m-1})$ a component of system responses to family (7) that reflects contribution of kernel $K_m$ (this component can be determined without any principal difficulties), then taking into account symmetry of $K_m$ with respect to all variables its identification is reduced to solving the linear N-dimensional Volterra equation of the first kind

$$V_m K_m \equiv \sum_{i_1 + \cdots + i_{m-1} = m} (-1)^{\sum_{k=1}^{\left[\frac{m-1}{2}\right]} i_{2k}} \frac{m!}{i_1! \ldots i_{m-1}!} V_{i_1, \ldots, i_{m-1}} K_m = f_m, \tag{10}$$

where

$$V_{i_1, \ldots, i_{m-1}} K_m = \underbrace{\int_{t-\omega_1}^{t} \cdots \int}_{i_1 \text{ раз}} \cdots \underbrace{\int_{t-\omega_1-\cdots-\omega_{m-1}}^{t-\omega_1-\cdots-\omega_{m-2}} \cdots \int}_{i_{m-1} \text{ раз}} K_m(s_1, \ldots, s_m) ds_1 \ldots ds_m,$$

$$f_m \equiv f_m(t, \omega_1, \ldots, \omega_{m-1}),$$

$$t, \omega_1, \ldots, \omega_{m-1} \in \Delta_m = \left\{ t, \omega_1, \ldots, \omega_{m-1} / 0 \leq \sum_{k=1}^{m-1} \omega_k \leq t \leq T \right\},$$

and symbol $[\ldots]$ in the limit superior of summation (10) means an integer part of the number. For $m = 2,3$ from (10) we have correspondingly

$$V_2 K_2 \equiv \int_{t-\omega_1}^{t} \int_{t-\omega_1}^{t} K_2(s_1, s_2) ds_1 ds_2 = f_2(t, \omega_1), \quad (11)$$

$$V_3 K_3 \equiv \overbrace{\int\int\int_{t-\omega_1}^{t}}^{t} K_3(s_1, s_2, s_3) ds_1 ds_2 ds_3 - 3 \overbrace{\int\int_{t-\omega_1}^{t}}^{t} \int_{t-\omega_1-\omega_2}^{t-\omega_1} K_3(s_1, s_2, s_3) ds_1 ds_2 ds_3 + \quad (12)$$

$$+ 3 \int_{t-\omega_1}^{t} \overbrace{\int\int_{t-\omega_1-\omega_2}^{t-\omega_1}}^{t-\omega_1} K_3(s_1, s_2, s_3) ds_1 ds_2 ds_3 - \overbrace{\int\int\int_{t-\omega_1-\omega_2}^{t-\omega_1}}^{t-\omega_1} K_3(s_1, s_2, s_3) ds_1 ds_2 ds_3.$$

In (Apartsyn, 1991; Apartsyn, 1992) the authors give the inversion formulas (11), (12):

$$K_2(t, t-\omega_1) = \frac{f''_{2_{t\omega_1}} + f''_{2_{\omega_1^2}}}{2}, \quad (13)$$

$$K_3(t, t-\omega_1, t-\omega_1-\omega_2) = \frac{f'''_{3_{t\omega_1^2}} - f'''_{3_{t\omega_1\omega_2}} + f'''_{3_{\omega_1\omega_1^2}} - f'''_{3_{\omega_1^2\omega_2}}}{12}. \quad (14)$$

For large values $m$ the inversion formulas (10) are presented in (Apartsyn, 1996; Apartsyn, 1999; Apartsyn, 2003). However, the applied significance of these formulas is not high due to instability of the operation of numerical differentiation of empiric functions. If to replace the variables $\omega_i = s_i - s_{i+1} \geq 0, i = \overline{1, m-1}, s_1 = t$, then, as is shown in (Apartsyn, 1996), the inversion formula (10) can be represented in a compact form:

$$K_m(s_1, \ldots, s_m) = \frac{(-1)^{\left[\frac{m}{2}\right]}}{m! \, 2^{m-2}} f^{(m)}_{m_{s_1,\ldots,s_m}}(s_1, s_1 - s_2, \ldots, s_{m-1} - s_m). \quad (15)$$

The case of vector input, which is important for applications, was studied in (Apartsyn, 1996; Apartsyn, Guseva, Solodusha, Hudyakov & Tairov, 1992; Apartsyn & Solodusha, 1997) discrete analogs of inversion formulas that are based on self-regularizing (Apartsyn, 1987) property of quadrature of middle rectangles were studied in (Apartsyn, Guseva, Solodusha, Hudyakov & Tairov, 1992; Apartsyn & Solodusha, 1997; Solodusha, 1992a), and in (Apartsyn, 1995b; Apartsyn, 1997; Sidorov, 1997; Sidorov, 1998) authors considered a non-stationary case when the kernels $K_m$ explicitly depend on $t$.

A series of studies (Apartsyn, Solodusha & Tairov, 1992; Apartsyn, Solodusha, Tairov & Khudyakov, 1994; Apartsyn, Solodusha & Tairov, 1997; Apartsyn & Tairov, 1996; Solodusha, 1992b; Solodusha, 1994) is devoted to application of the presented approach in modeling of heat exchange processes in the high-temperature plant at the Energy Systems Institute SB RAS. S.V. Solodusha (Solodusha, 2012a) has developed a software package for modeling nonlinear dynam-ics of heat exchange. The software was used to process real experimental data. Description of the software is presented in the next section.

It should be noted that along with obvious advantages which follow from availability of explicit inversion formulas of respective multidimensional integral Volterra equations of the first kind the developed method of identifying the Volterra kernels has a flaw related to quite severe conditions for the existence of solutions to these equations in the required classes of functions (Apartsyn, 1995c; Apartsyn, 1996). Therefore, the authors of (Apartsyn, 2005) note, that in order to model a response of system $y(t)$ to an input disturbance $x(t)$ the knowledge of $K_m$, is in general redundant. It suffices to calculate multi-dimensional convolutions in (3). If for their approximate calculation to apply the product integration method (Linz, 1985), according to which at a chosen step of net $h$ in the one-dimensional case

$$\int_0^{ih} K_1(s)x(t-s)ds \approx \sum_{j=1}^{i} x(t_{i-j+\frac{1}{2}}) \int_{(j-1)h}^{jh} K_1(s)ds, \quad (16)$$

$t_{i-\frac{1}{2}} = \left(i - \frac{1}{2}\right)h, i = \overline{1,n}, nh = T$, then, instead of identification of $K_m$ themselves, it is sufficient of identify elementary integrals of $K_m$. In (Apartsyn & Spiryaev, 2005) this approach is implemented for the case $K_2(s_1, s_2) = \psi(s_1)\psi(s_2)$, $K_3(s_1, s_2, s_3) = \varphi(s_1)\varphi(s_2)\varphi(s_3)$, $\psi(t), \varphi(t) \in C_{[0,T]}$, and in (Spiryaev, 2006) test signals (8), (9) were used to obtain explicit inversion formulas of corresponding systems of linear algebraic equations for a general case where $K_2, K_3$ are random continuous symmetric functions.

Let us dwell on the important problem of increasing the accuracy of modeling the nonlinear dynamics by Volterra polynomials through optimization of test signal amplitude. A "reference" mathematical model was introduced in (Apartsyn, 2001; Apartsyn & Solodusha, 2004)

$$y_{ref}^N(t) = \sum_{m=1}^{N} \frac{1}{m!} \left(\int_0^t x(s)ds\right)^m, \quad (17)$$

which is an $N$- section of the series

$$y_{ref}(t) = \sum_{m=1}^{\infty} \frac{1}{m!} \left(\int_0^t x(s)ds\right)^m = e^{\int_0^t x(s)ds} - 1. \quad (18)$$

If in (17) $N \geq 3$, an error in modeling (17) with quadratic Volterra polynomial

$$y_{sq}(t) = P_2(x(t)) \equiv \int_0^t K_1(s)x(t-s)ds + \int_0^t \int_0^t K_2(s_1, s_2)x(t-s_1)x(t-s_2)ds_1 ds_2, \quad (19)$$

is non-zero. Moreover, it depends on both the feasible set of input signals and the chosen amplitudes $\alpha_i, i = \overline{1,2}$, of test signals (8).

At $N \geq 4$ in a general case an error in modeling (17) is non-zero when modeling is performed with the cubic Volterra polynomial

$$y_{cub}(t) = P_3(x(t)) \equiv P_2(x(t)) + \int_0^t \int_0^t \int_0^t K_3(s_1, s_2, s_3) \prod_{i=1}^3 x(t - s_i) \, ds_i, \qquad (20)$$

and the choice of $\alpha_i, i = \overline{1,3}$ in (9) is significant again.

In (Apartsyn, 1995c) the authors established that for solvability of the integral equations (11), (12) in the class of continuous symmetric kernels $K_2$ and $K_3$ it is necessary to meet the conditions

$$\alpha_1 + \alpha_2 = 0 \qquad (21)$$

for (8) and

$$\alpha_1 + \alpha_2 + \alpha_3 = 0 \qquad (22)$$

for (9).

Since equalities (21), (22) maintain arbitrariness while choosing concrete real non-zero amplitudes of test signals (8), (9), it is natural to formulate some extreme problems of minimizing the errors between $y_{ref}(t)$ and models $y_{sq}(t)$ and $y_{cub}(t)$.

First consider the quadratic model (19). Choose an elementary set of signals

$$X(B, T) = \{x^\beta(t) = \beta I(t); \beta \in (0, B], t \in [0, T]\}, \qquad (23)$$

to be admissible and the absolute value of residual between (17) and $y_{sq}(t)$ at $t = T$:

$$\left\| [\mathcal{N}_{N,sq}(x^\beta(t))] \right\|_{t=T} \triangleq \left| y_{et}^{N,\beta}(T) - y_{sq}^{\alpha,\beta}(T) \right| \qquad (24)$$

to be an objective function.

In (24) $y_{sq}^{\alpha,\beta}(t)$ is a response to signal $x^\beta(t)$ of model (19) in which kernels $K_1$ and $K_2$ were identified using a family of test disturbances (8) with $\alpha_1 = \alpha$ and $\alpha_2 = -\alpha$, whereas $y_{ref}^{N,\beta}(t)$ is a response of the reference model (17) to $x^\beta(t)$.

Choice of criterion (24) is dictated by the fact that in many applications it is important to know how to estimate the value of response of a dynamic object at the end of the transient process.

In fact residual $\mathcal{N}_{N,sq}(x^\beta(T))$ at fixed $N, T, B$ is the function of parameters $\alpha$ and $\beta$: $\mathcal{N}_{N,sq} = \mathcal{N}_{N,sq}(\alpha, \beta)$. Thus, we arrive at the following extreme problem: find

$$\alpha_N^* = \arg \min_{\alpha \in (0,B]} \left\{ \max_{\beta \in (0,B]} |\mathcal{N}_{N,sq}(\alpha, \beta)| \right\}, N \geq 3. \qquad (25)$$

At $N = 3$ problem (25) can be solved in an explicit form. Indeed, it is easy to check that

$$\mathcal{N}_{3,sq}(\alpha, \beta) = \frac{(\beta^3 - \alpha^2 \beta)T^3}{6},$$

from which $\beta_{max} = \left\{ \frac{B}{2}, B \right\}$; $\alpha_3^* = \frac{\sqrt{3}}{2} B \approx 0.866B$ and

$$\mathcal{N}_{3,sq}(\alpha_3^*, \beta_{max}) = \frac{B^3 T^3}{24}. \tag{26}$$

Although for $N > 3$ the extreme problem (25) can be solved only numerically, the calculations show that with increasing $N$ the values of $\alpha_N^*$ get stabilized fast at $N \geq 6$ $\alpha_N^* \approx 0.878B$ as well. If to consider similar extreme problem for the quadratic model constructed on the basis of the product integration method:

$$y_{sq,pi}(ih) = \sum_{j=1}^{i} m_j x\left(t_{i-j+1/2}\right) + \sum_{j=1}^{i}\sum_{k=1}^{i} l_{jk} x\left(t_{i-j+1/2}\right) x\left(t_{i-k+1/2}\right), \tag{27}$$

where

$$m_j = \int_{(j-1)h}^{jh} K_1(s)ds; \quad l_{jk} = \int_{(j-1)h}^{jh}\int_{(k-1)h}^{kh} K_2(s_1, s_2) ds_1 ds_2; \quad i = \overline{1,n}, nh = T, \tag{28}$$

then by virtue of the absence of constraint (21) the accuracy of (27) as compared to (19) can be increased by choosing the amplitudes $\alpha_1$ and $\alpha_2$ from the section $(0, B]$. Indeed, the calculation shows (Apartsyn & Spirayev, 2005) that

$$\mathcal{N}_{3,sq,pi}(\alpha_1, \alpha_2, \beta) = (\beta^3 - (\alpha_1 + \alpha_2)\beta^2 + \alpha_1\alpha_2\beta)\frac{T^3}{6}, \tag{29}$$

therefore, the problem of search

$$(\alpha_{3,1}^*, \alpha_{3,2}^*) = \arg\min_{\alpha_1,\alpha_2 \in (0,B]} \left\{\max_{\beta \in (0,B]} |\beta^3 - (\alpha_1 + \alpha_2)\beta^2 + \alpha_1\alpha_2\beta|\right\}$$

has an accurate solution

$$\alpha_{3,1}^* = (2\sqrt{3} - 3)B \approx 0.464; \quad \alpha_{3,2}^* = 2\alpha_{3,1}^* \approx 0.928;$$

$$\beta_{max} = \{(3\sqrt{3} - 5)B, (\sqrt{3} - 1)B, B\},$$

here

$$\mathcal{N}_{3,sq,pi}(\alpha_{3,1}^*, \alpha_{3,2}^*, \beta_{max}) \approx 0.0064 B^3 T^3,$$

which is one order less than (26).

If to expand the admissible set $X(B,T)$, by replacing condition $\beta \in (0, B]$ by $\beta \in [-B, B]$, then from the symmetry of (29) with respect to $\alpha_1, \alpha_2$ it follows that $\alpha_{3,1}^* = -\alpha_{3,2}^* = \alpha_3^* = \frac{\sqrt{3}}{2}B$ and models (19), (27) are comparable in accuracy.

In the end, the expansion of $X(B,T)$ by including signals of type (9)

$$x_{\omega_1,\omega_2}^{\beta}(t) = \beta(I(t) - 2I(t - \omega_1) + I(t - \omega_1 - \omega_2))$$

results, as is shown in (Apartsyn & Solodusha, 2004), in the following expression for $\mathcal{N}_{3,sq}(\alpha, \beta, \omega_1, \omega_2)$:

$$\mathcal{N}_{3,sq}(\alpha, \beta, \omega_1, \omega_2) = \frac{(\beta^3 - \beta\alpha^2)(\omega_1 - \omega_2)^3}{6} - \beta\alpha^2 \omega_1 \omega_2^2,$$

and numerical solution to the problem of search

$$\alpha_3^* = arg \min_{\alpha \in (0,B]} \left\{ \max_{\substack{0 \leq \omega_1, \omega_2 \leq T \\ \beta \in [-B,B]}} |\mathcal{N}_{3,sq}(\alpha, \beta, \omega_1, \omega_2)| \right\}$$

is as follows:

$$\alpha_3^* \approx 0.732B; \; \omega_{1_{max}} = 0.366T; \; \omega_{2_{max}} = 0.634T; \; \beta_{max} = B.$$

Proceeding to the case of cubic models, assume in (17) $N = 4$ and consider the problem of search for a pair, subject to (22),

$$(\alpha_{4,1}^*, \alpha_{4,2}^*) = arg \min_{\alpha_1, \alpha_2 \in (0,B]} \left\{ \max_{\beta \in (0,B]} |\mathcal{N}_{4,cub}(\alpha_1, \alpha_2, \beta)| \right\}, \tag{30}$$

where, as is established in (Apartsyn & Solodusha, 2004), a residual between $y_{ref}^{4,\beta}(T)$ and response (20) $y_{cub}^{\alpha_1, \alpha_2, \beta}(T)$ to signal $x^\beta(t) = \beta I(t)$ has the form

$$\mathcal{N}_{4,cub}(\alpha_1, \alpha_2, \beta) = \{\beta^4 + [\alpha_1 \alpha_2 - (\alpha_1^2 + \alpha_1 \alpha_2 + \alpha_2^2)]\beta^2 + \alpha_1 \alpha_2 (\alpha_1 + \alpha_2)\beta\} \frac{T^4}{24}. \tag{31}$$

The numerical solution to (30), (31), which was obtained with the MAPLE system is:

$$\alpha_{4,1}^* \approx 0.475B, \alpha_{4,2}^* \approx 0.885B, \beta_{max} = B,$$

$$\mathcal{N}_{4,cub}(\alpha_{4,1}^*, \alpha_{4,2}^*, \beta_{max}) \approx 0.00037 B^4 T^4.$$

Emphasize the significance of holding the equality (22) for model (20). If it does not hold, then, as is shown on a concrete example in (Apartsyn & Solodusha, 2004), the cubic model can give even much worse results than the quadratic model (19). If $\beta \in (0, B]$, the absence of the need to meet condition (22) again gives advantage to the cubic pi-model, for which the residual (Apartsyn & Spiryaev, 2005) is

$$\mathcal{N}_{4,cub,pi}(\alpha_1, \alpha_2, \alpha_3, \beta) =$$

$$= \{\beta^4 - (\alpha_1 + \alpha_2 + \alpha_3)\beta^3 + (\alpha_1 \alpha_2 + \alpha_1 \alpha_3 + \alpha_2 \alpha_3)\beta^2 - \alpha_1 \alpha_2 \alpha_3 \beta\} \frac{T^4}{24},$$

and numerical values of optimal amplitudes are

$$\alpha_{4,1}^* \approx 0.283B, \; \alpha_{4,2}^* \approx 0.677B, \; \alpha_{4,3}^* \approx 0.960B.$$

Finally, to sum up this section note that although the problem of optimizing the amplitudes of test signals of type (7) was considered above as applied to the reference model (17), the obtained recommendations are of a more universal character than it is shown in (Apartsyn & Solodusha, 2004) on the example of the applied problem of modeling the nonlinear dynamic processes of heat exchange.

## SOFTWARE FOR MODELING NONLINEAR DYNAMICS

A software package based on the reference model was created to construct and test integral models of nonlinear dynamics of heat exchange. A description of heat exchange process that occurs in the component of heat exchanger with independent heat supply (Tairov, 1989) was used as a reference. According to (Tairov, 1989) deviation of enthalpy at the output $\Delta i(t)$ at random laws of disturbances of liquid flow rate $\Delta D(t)$ and heat supply $\Delta Q(t)$ is defined by the relationship:

$$\Delta i(t) = \frac{\lambda_1 \lambda_2}{\lambda_2 - \lambda_1} \int_0^t \left( \Delta Q(s) - \frac{Q_0}{D_0} \Delta D(s) \right) \left( e^{-\lambda_1 \int_s^t D(s_1) ds_1} - e^{-\lambda_2 \int_s^t D(s_1) ds_1} \right) ds. \qquad (32)$$

In (32) $t$ – time, $\lambda_1$ and $\lambda_2$ – some constants, indices '0' are used to denote parameters of an initial stationary condition, $\Delta$ – increment to a corresponding parameter of the initial stationary condition, for example $D(t) = D_0 + \Delta D(t)$.

A software package was created in the object-oriented programming environment Borland C++ Builder and is based on the function-module principle. The software package includes the following modules:

- Modules of constructing integral models in the form of linear, quadratic and cubic Volterra polynomials for scalar input signals;
- Modules of constructing integral models in the form of linear and quadratic Volterra polynomials for vector input signals;
- Module of calculating a controlled input action that provides a reliable (specified) response of the integral model (such problem arises in relation to the problems of automatic regulation of technical objects).

The software consists of blocks intended for adjustment of input parameters, identification and modeling.

In a dialog mode user can change the reference model parameters, and specify discretization step, length of time section $[0, T]$, and amplitude of test input signals.

The identification block calculates net analogues of the reference model responses to sets of test input signals. The obtained data are used for approximation of transient characteristics of dynamic object. The procedures for identification are based on difference analogs of inversion formulas, which provides fast on-line operation.

In the modeling block the integral models are used to calculate the dynamic processes for input signals of arbitrary form. The results are compared with an output signal of the reference model (32). The procedures for calculation of responses of integral models are based on the methods of middle rectangles and product integration. Output values of the reference model are calculated with the help of the trapezium method. The software implements the input functions and digitiza-tion of input actions, as well as graphical representation of computational experiment results. User can study integral models by choosing input disturbances from a database or by setting signals with the help of manipulator "mouse". All the output data are stored in respective files on a disk and can be used for detailed familiarization and analysis.

This software was applied in (Solodusha, Spiryaev & Scherbinin, 2006) to model nonlinear dynamics of heat exchange by the cubic Volterra polynomial in a scalar case. The tactics of choosing the test input signals in (Solodusha, Spiryaev & Scherbinin, 2006) differed from that previously used in (Apartsyn, Solodusha & Sidorov, 1998) and the model construction itself was based on the product integration method. In this paper we will demonstrate the software only for the case of vector input $(x(t) = (\Delta D(t), \Delta Q(t)))$ of a quadratic model.

Figures 1(a) and 2(a) show arbitrary input disturbances $\Delta D(t)$ and $\Delta Q(t)$ in percentage of their stationary values $D_0 = 0.16 \ (kg/s), Q_0 = 100$ (kW). The results of modeling $\Delta i(t)$(kJ/kg) on the basis of reference model and integral models for specified input signals are presented in Figures 1 (b) and 2 (b). Construction of integral models is based on the algorithms described in (Apartsyn, Tairov, Solodusha, Khudyakov, 1994) and (Solodusha, Spiryaev & Scherbinin, 2007), respectively. Both integral models employed identical amplitudes of test signals, that are equal to $25\%D_0$ and $25\%Q_0$.

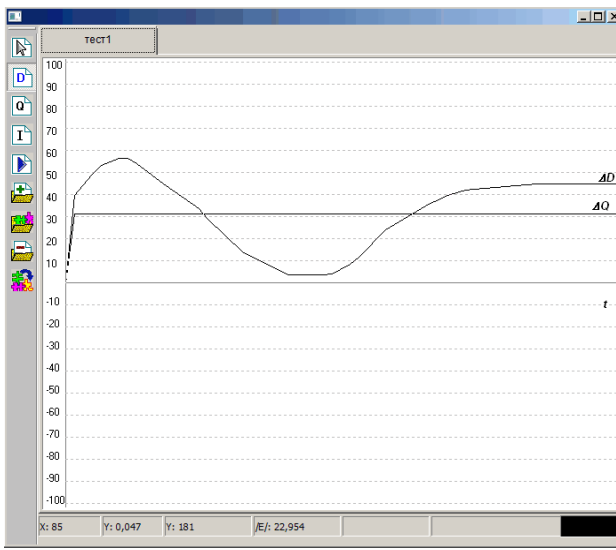 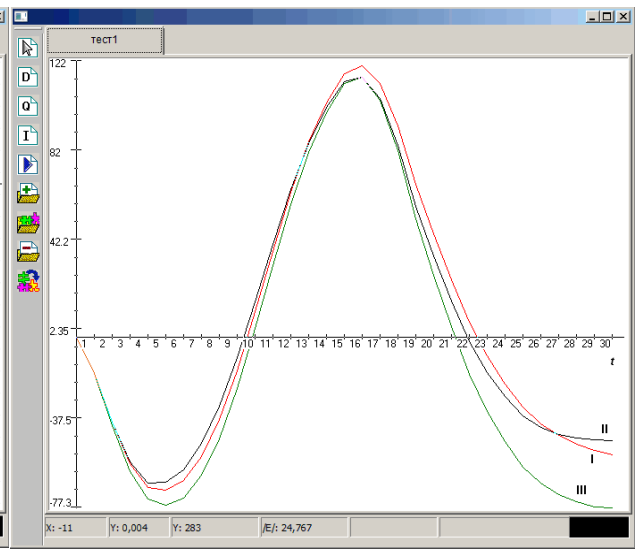

*Figure. 1(a). Curves of input disturbances.*     *Figure. 1(b). Modeling results based on the reference model (I) and integral models (II, III).*

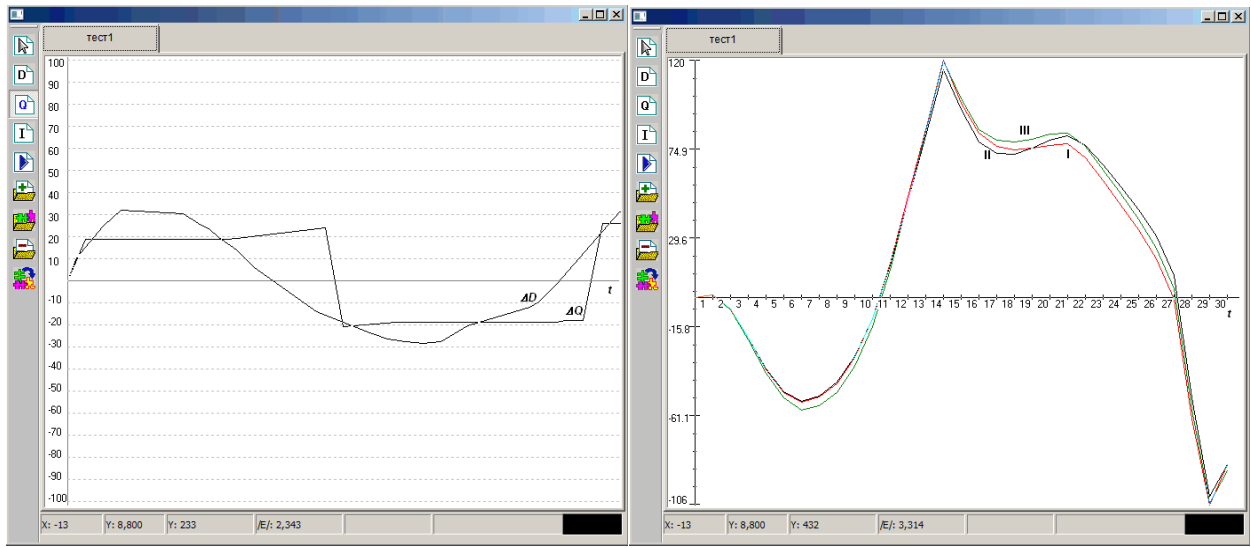

*Figure. 2(a). Curves of input disturbances.*   *Figure. 2(b).Modeling results based on the reference model (I) and integral models (II, III).*

In the numerical experiments $h = 1(\text{sec}), t \in [0,30]$. Such a choice of parameters is related to real data obtained in the course of experiment conducted on the high-temperature plant at ESI SB RAS. Results of the calculations show that both models give an acceptable accuracy of modeling the heat exchange dynamics.

## POLYNOMIAL VOLTERRA EQUATIONS OF THE FIRST KIND AND LAMBERT FUNCTION

Until now the problem of identifying the Volterra kernels $K_m, m = \overline{1,N}$ in $(2), (3)$ has been the focal point of the research.

Assume now that this problem has already been solved with some approach. Then, the stage of mathematical modeling to be considered next can be a typical problem of automatic control of a dynamic object in the absence of feedback, i.e. finding an input signal $x(t)$ such that it corresponds to a specified (desired) response $y(t)$. At known $K_m, m = \overline{1,N}$, and $y(t)$ (2), (3) are non-linear, integral equations with respect to $x(t)$. Emphasize that although equations (2) and (3) include multi-dimensional integrals of Volterra type the sought function is the function of only one variable – time $t$. It seems surprising that with abundant publications on the methods for identification of Volterra polynomials such integral equations have not been considered in the literature until recently. Moreover, there is still no fixed name for them. In (Apartsyn, 2000; Apartsyn, 2004a; Apartsyn, 2004b; Apartsyn, 2005; Apartsyn, 2007) (2) and (3) are interpreted as multilinear (in particular at $N = 2,3$ – bilinear and threelinear, respectively) Volterra equations of the first kind. In (Belbas & Bulka, 2011) the authors suggest the term multiple nonlinear, and in (Apartsyn, 2010; Apartsyn, 2011; Apartsyn, 2012) the authors use the name "polynomial" which seems the most appropriate since it reflects the fact that to solve (2) and (3) implies to find the "roots" (or the only root in the space of real continuous functions) of the Volterra polynomial.

It is well known (Apartsyn, 1999; Apartsyn, 2003), that at $N = 1$ in (2) the linear equation

$$P_1(x(t)) \equiv \int_0^t K_1(t,s)x(s)ds = y(t), \quad t \in [0,T], \tag{33}$$

on the assumption that $K_1(t,t) \neq 0 \;\forall t \in [0,T]; K_{1_t}'(t,s) \in C_\Delta$, $\Delta = \{t,s/0 \leq s \leq t \leq T\}$; $y(t) \in \dot{C}_{[0,T]}^{(1)}$ $(y(0) = 0)$ has the unique continuous solution $x(t) \in C_{[0,T]}$ at any $T < \infty$.

It appears that the main specific feature of the polynomial equation (2) at $N > 1$ lies in the fact that the existence of a real continuous solution to (2) is guaranteed at sufficiently small $T > 0$, i.e. this solution is of a local character.

For understanding the specific features of the polynomial equation (2) assume $N = 2$ and consider an elementary quadratic Volterra equation of the first kind with constant kernels

$$K_1(t,s) \equiv 1, \quad K_2(t,s_1,s_2) \equiv \lambda \neq 0,$$

so that (2) has the form

$$P_2(x(t)) \equiv \int_0^t x(s)ds + \lambda \left( \int_0^t x(s)ds \right)^2 = y(t), \quad t \in [0,T]. \tag{34}$$

Let the right-hand part of (34), along with conditions presented above for linear case, additionally satisfy the inequality

$$4\lambda y(t) > -1 \;\forall t \in [0,T]. \tag{35}$$

It is easy to show that condition (35) is necessary and sufficient for the existence and uniqueness of the real continuous solution $x^*(t) \in C_{[0,T]}$ at any $0 < T < \infty$, in this case the following inversion formula (34) is true:

$$x^*(t) = P_2^{-1}(y(t)) = \frac{y'(t)}{\sqrt{1 + 4\lambda y(t)}}, t \in [0,T]. \tag{36}$$

It is clear that if $y(t)$ has a constant sign on $[0,T]$ and $sign\,(y(t)) = sign\,(\lambda)$, then (35) a fortiori holds at any $0 < T < \infty$. However, if signs $y(t)$ and $\lambda$ are opposite, or $y(t)$ is alternating on $[0,T]$, then a sufficient condition for the existence of the real continuous solution to (2) at any $0 < T < \infty$ is inequality

$$|y(t)| < \frac{1}{4|\lambda|}, \quad \forall t \in [0,T]. \tag{37}$$

Since $y(0) = 0$, then in the neighborhood of point $t = 0$

$$y(t) = ty'(\xi), 0 \leq \xi \leq t, \tag{38}$$

and it means that (37) follows from inequality

$$tF(t) < \frac{1}{4|\lambda|}, \quad \forall t \in [0, T], \tag{39}$$

in which

$$F(t) = \max_{0 \leq \xi \leq t} |y'(\xi)|. \tag{40}$$

However, (39) cannot hold at any $0 < T < \infty$, since $F(t)$ is continuous and monotonously non-decreasing and, therefore, the function $tF(t)$ is continuous and strictly increasing, equal to zero at $t = 0$, therefore, a fortiori there will be $t = T^*$, such that

$$T^* F(T^*) = \frac{1}{4|\lambda|}. \tag{41}$$

Thus, condition (37) is met for sure only at $T < T^*$. The issue of estimating the maximum possible region of definition of the real continuous solution to (2) at set $K_m$ and $y(t)$ in (2), (3) turns out to be a key one in the theory of polynomial Volterra equations of the first kind.

The Lambert function plays an important part for further considerations. Wide application of the Lambert function in the mathematical studies became possible owning to the developers of the MAPLE system (Corless, Gonnet, Hare, Jeffery & Knuth, 1996; Corless, Gonnet, Hare & Jeffrey, 1993; Corless, Jeffrey & Knuth, 1997). Below the Lambert function properties that are most important for this paper are presented in brief.

Consider the function of real argument $z$

$$y(z) = ze^z, z \in (-\infty, \infty). \tag{42}$$

Solution to the equation

$$ze^z = y,$$

is called the Lambert function, i.e. $W(y)$ – the function inverse to function (42). The Lambert function has two real branches – the main, which is defined for $y \in \left[-\frac{1}{e}, \infty\right]$ and is analytic at zero (it is denoted by $W_0(y)$ or $W(y)$), and the second branch defined for $y \in \left[-\frac{1}{e}, 0\right]$ and denoted by $W_{-1}(y)$. At $y = -\frac{1}{e}$ $W\left(-\frac{1}{e}\right) = W_{-1}\left(-\frac{1}{e}\right) = -1$, and at $y = 0$ $W(0) = 0, W_{-1}(0) = -\infty$. Both branches at point $y = -\frac{1}{e}$ have a vertical tangent line. Directly from the definition it follows that

$$W'(y) = \frac{W(y)}{(1 + W(y))y}, y > -\frac{1}{e}. \tag{43}$$

The properties of the Lambert function in complex domain as well as its applications in many areas of the applied mathematics are described in detail in (Corless, Gonnet, Hare, Jeffery & Knuth, 1996; Corless, Gonnet, Hare & Jeffrey, 1993; Corless, Jeffrey & Knuth, 1997; Apartsyn, 2008).

Note, that since $e^z = \sum_{v=0}^{\infty} \frac{z^v}{v!}$, the function $W(y)$ is a solution to the problem of inversion of the infinite series

$$\sum_{v=1}^{\infty} \frac{z^v}{(v-1)!} = y,$$

Besides, the representation (Corless, Jeffrey & Knuth, 1997)

$$W(y) = \sum_{k=0}^{\infty} \frac{(-k)^{k-1}}{k!} y^k = y - y^2 + \frac{3}{2}y^3 - \frac{8}{3}y^4 + \frac{125}{24}y^5 + O(y^6) \qquad (44)$$

is true.

The previously introduced reference mathematical model $y_{et}(t)$ (see (18)) corresponds to the Volterra kernels

$$K_m = \frac{1}{m!}, m = 1,2, ....$$

It is convenient to denote

$$\Theta(t) = \int_0^t x(s)ds, t \in [0,T], \qquad (45)$$

then, according to (18),

$$y_{et}(t) + 1 = e^{\Theta(t)},$$

so that

$$\Theta(t) = \ln(1 + y_{et}(t)),$$

from which at specified $y_{et}(t)$ the real continuous input signal $x_\infty^*(t)$ is determined by the following inversion formula of the series (18):

$$x_\infty^*(t) = \Theta'(t) = \frac{y_{et}'(t)}{1 + y_{et}(t)}, t \in [0,T]. \qquad (46)$$

In (46) $y_{et}(t) \in C_{[0,T]}^{(1)}, y_{et}(0) = 0, y_{et}(t) > -1 \; \forall t \in [0,T]$.

Thus, if $y_{et}(t) = t$, then, according to (46),

$$x_\infty^*(t) = \frac{1}{1+t}, t \in [0,T],$$

and at $y_{et}(t) = -t$

$$x_\infty^*(t) = \frac{1}{t-1} \notin C_{[0,T]},$$

if $T \geq 1$, since the condition $y_{et}(t) > -1 \; \forall t \in [0,T]$ is not met.

Let now in (1)

$$K_m = \frac{1}{(m-1)!}, m = 1,2, ..., \qquad (47)$$

so that

$$y(t) = \sum_{m=1}^{\infty} \frac{\left(\int_0^t x(s)ds\right)^m}{(m-1)!} = \sum_{m=1}^{\infty} \frac{\Theta^m(t)}{(m-1)!} = \Theta(t)e^{\Theta(t)}. \tag{48}$$

It follows from (45) that $\Theta(0) = 0$, therefore according to the definition of the main branch of the Lambert function

$$\Theta(t) = W_0(y(t)), \tag{49}$$

from which taking into account (43)

$$x^*(t) = \frac{W_0(y(t))y'(t)}{\left(1 + W_0(y(t))\right)y(t)}, t \in [0,T]. \tag{50}$$

In (50) $y(t) \in C_{[0,T]}^{(1)}$, $y(0) = 0$, $y(t) > -\frac{1}{e}$ $\forall t \in [0,T]$, and it follows from (44) that $x_\infty^*(0) = y'(0)$. If, for example $y(t) = te^t$, then according to (50) $x_\infty^*(t) = 1$, but if $y(t) = -t$, then at $T \geq -\frac{1}{e}$ the function

$$x_\infty^*(t) = \frac{W_0(-t)}{(1 + W_0(-t))t}$$

is not continuous on $[0,T]$ due to discontinuity of the second kind at point $t = -\frac{1}{e}$ $\left(W_0\left(-\frac{1}{e}\right) = -1\right)$.

The Lambert function naturally occurs when studying the polynomial Volterra equations of the first kind (2). Considering the quadratic equation (34) we have found out that the critical $T = T^*$ satisfies equation (41). If $y(t) = e^t - 1$, then by virtue of (40) $F(t) = e^t$ and (41) gives $T^*e^{T^*} = \frac{1}{4|\lambda|}$, from which

$$T^* = W_0\left(\frac{1}{4|\lambda|}\right).$$

Consider the quadratic equation of Volterra of the first kind which is more general than (34) by assuming in (2) at $N = 2$

$$K_1(t,s) = 1 - L_1(t-s), L_1 > 0; K_2(t,s_1,s_2) \equiv -\lambda, \lambda > 0,$$

so that (2) has the form

$$P_2(x(t)) \equiv \int_0^t (1 - L_1(t-s))x(s)ds - \lambda\left(\int_0^t x(s)ds\right)^2 = y(t), \ t \in [0,T]. \tag{51}$$

Let $y(t) = \mathcal{F}t, \mathcal{F} > 0$. By differentiating (51) proceed to the equivalent equation of the second kind

$$x(t) - L_1\int_0^t x(s)ds - 2\lambda x(t)\int_0^t x(s)ds = \mathcal{F}, \tag{52}$$

from which, using replacement of (45), arrive at the Cauchy problem for the nonlinear differential equation

$$\dot{\Theta}(t) = \frac{\mathcal{F} + L_1\Theta(t)}{1 - 2\lambda\Theta(t)}, \Theta(0) = 0, t \in [0, T]. \tag{53}$$

Solution to (53) $\Theta_2^*(t)$ is expressed through the main real branch of the Lambert function:

$$\Theta_2^*(t) = -\frac{L_1 + 2\lambda\mathcal{F}}{2\lambda L_1} W_0\left(-\frac{2\lambda\mathcal{F}}{L_1 + 2\lambda\mathcal{F}} \exp\frac{L_1 t - 2\lambda\mathcal{F}}{L_1 + 2\lambda\mathcal{F}}\right) - \frac{\mathcal{F}}{L_1},$$

from which, considering (43), obtain the following inversion formula (51) at $y(t) = \mathcal{F}t$:

$$x_2^*(t) = P_2^{-1}(\mathcal{F}t) = -\frac{L_1}{2\lambda}\frac{W_0\left(-\frac{2\lambda\mathcal{F}}{L_1 + 2\lambda\mathcal{F}} \exp\frac{L_1 t - 2\lambda\mathcal{F}}{L_1 + 2\lambda\mathcal{F}}\right)}{1 + W_0\left(-\frac{2\lambda\mathcal{F}}{L_1 + 2\lambda\mathcal{F}} \exp\frac{L_1 t - 2\lambda\mathcal{F}}{L_1 + 2\lambda\mathcal{F}}\right)}, \tag{54}$$

and equating the argument $W_0$ to $-\frac{1}{e}$ gives the upper estimate for $T$ in (52), (53):

$$T < T^* = \frac{L_1 + 2\lambda\mathcal{F}}{L_1^2}\ln\left(1 + \frac{L_1}{2\lambda\mathcal{F}}\right) - \frac{1}{L_1}. \tag{55}$$

Let us consider the case of random $N \geq 1$ in (2).

As is mentioned above at $N = 1$ the linear integral equation (33) on the natural assumptions about the kernel $K(t, s)$ and right-hand part $y(t)$ has the unique continuous solution $x_1^*(t)$ at any $0 < T < \infty$. Moreover, it is stable to errors in the right-hand part which are measured in the norm $\dot{C}_{[0,T]}^{(1)}$, which follows from the estimation in (Apartsyn, 1999; Apartsyn, 2003)

$$\|P_1^{-1}(y(t))\|_{\dot{C}_{[0,T]}^{(1)} \to C_{[0,T]}} \leq \frac{1}{k}e^{L_1 k^{-1}T}, \tag{56}$$

where

$$k = \min_{t \in [0,T]}|K_1(t, t)| > 0, \tag{57}$$

$$L_1 = \max_{0 \leq s \leq t \leq T}|K'_{1_t}(t, s)| \geq 0. \tag{58}$$

This means that equation (33) is well-posed on the pair $\left(C_{[0,T]}, \dot{C}_{[0,T]}^{(1)}\right)$ according to Hadamard.

Suppose additionally that the kernels $K_m$, $m = \overline{2, N}$, in (2) are not only symmetrical with respect to $s_1, \ldots, s_m$ and continuous with respect to a set of variables, but also continuously differentiable with respect to $t$. We are interested in the real continuous solution to (2) on the segment $[0, T]$. It is obvious that at small $T$

$$\left|\sum_{m=2}^{N}\int_0^t\ldots\int_0^t K_m(t, s_1, \ldots, s_m)\prod_{i=1}^m x(s_i)ds_i\right| = O(T^2),$$

therefore, equation (2) can be interpreted as a linear equation of type (33) with disturbed right-hand part $\tilde{y}(t)$, for which

$$\|y(t) - \tilde{y}(t)\|_{C^{(1)}_{[0,T]}} = O(T),$$

and which means that smallness of $T$ provides by virtue of stability of the solution to the linear equation the existence and uniqueness of the real continuous solution $x_N^*(t)$ on $[0,T]$ to equation (2). The principal point is the point about the estimation of value $T$, which guarantees this property. The universal technique that makes it possible to obtain such an estimate is the contraction mapping principle. However, as is shown in (Apartsyn, 2005) and (Apartsyn, 2007) it gives a too pessimistic upper estimate for $T$. Another approach which gives a guaranteed lower estimate of $T$ is based on introduction of special majorant integral equations and respective Cauchy problems.

Denote

$$L_m(t) = \max_{0 \leq s_1,\ldots,s_m \leq \xi \leq t} |K'_{m_t}(\xi, s_1, \ldots, s_m)| \geq 0, m = \overline{1,N}, t \in [0,T]; \quad (59)$$

$$M_m(t) = \max_{0 \leq s_2,\ldots,s_m \leq \xi \leq t} |K_m(\xi, \xi, s_2, \ldots, s_m)| \geq 0, m = \overline{2,N}, M_N > 0, t \in [0,T]. \quad (60)$$

Since (2) is equivalent, taking into account the condition $K_1(t,t) = 1 \,\forall t \in [0,T]$ (which does not decrease generality) and the symmetry of $K_m$ with respect to $s_1, \ldots, s_m$, to the integral Volterra equation of the second kind

$$x(t) + x(t) \sum_{m=2}^{N} m \int_0^t \ldots \int_0^t K_m(t,t,s_2,\ldots,s_m) \prod_{j=2}^{m} x(s_j) ds_j +$$

$$+ \sum_{m=1}^{N} \int_0^t \ldots \int_0^t K'_{m_t}(t,s_1,\ldots,s_m) \prod_{j=1}^{m} x(s_j) ds_j = y'(t), \quad t \in [0,T], \quad (61)$$

then $x_N^*(t)$ by virtue of (59), (60) satisfies the inequality

$$|x_N^*(t)| \leq F(t) + |x_N^*(t)| \sum_{m=2}^{N} m M_m(t) \left( \int_0^t |x_N^*(s)| ds \right)^{m-1} +$$

$$+ \sum_{m=1}^{N} L_m(t) \left( \int_0^t |x_N^*(s)| ds \right)^m, \quad t \in [0,T], \quad (62)$$

where $F(t)$ is determined in (40).

By replacing in (62) the sign $\leq$ by $=$, consider the integral equation

$$\psi(t) = F(t) + \psi(t) \sum_{m=2}^{N} m M_m(t) \left( \int_0^t \psi(s) ds \right)^{m-1} + \sum_{m=1}^{N} L_m(t) \left( \int_0^t \psi(s) ds \right)^m \quad (63)$$

and (taking into account the replacement of type (45)) the corresponding Cauchy problem

$$\dot{\Theta}(t) = G(\Theta(t)) \equiv \frac{F(t) + \sum_{m=1}^{N} L_m(t) \Theta^m(t)}{1 - \sum_{m=2}^{N} m M_m(t) \Theta^{m-1}(t)}, \Theta(0) = 0, t \in [0,T]. \quad (64)$$

In (Apartsyn, 2011) the authors show that (64) at sufficiently small $T$ has the unique solution $\Theta_N^*(t)$, and, hence, (63) has the unique solution $\psi_N^*(t) = \dot{\Theta}_N^*(t) \in C_{[0,T]}$. Besides from the results (Apartsyn, 2008) (lemma 1) it follows that the inequality

$$|x_N^*(t)| \leq \psi_N^*(t), t \in [0, T). \tag{65}$$

holds.

By virtue of (65) equation (63) and the Cauchy problem (64) can be considered majorant. For the case $F(t) = const > 0$, $L_m(t) = const \geq 0$, $M_m(t) = const > 0$ analogous majorant problems are introduced in (Apartsyn, 2004b; Apartsyn, 2005; Apartsyn, 2007).

In order to determine the maximum $T$ it suffices to equate the denominator in (64), in which $\Theta(t) \equiv \Theta_N^*(t)$, to zero.

The nonlinear functional equation with respect to $t$

$$\sum_{m=2}^{N} m M_m(t) [\Theta_N^*(t)]^{m-1} = 1 \tag{66}$$

has a unique solution, since by virtue of nondecrease $M_m(t)$ and strict increase $\Theta_N^*(t)$ the left-hand part of (66) is a continuous strictly increasing function $t$, equal to zero at $t = 0$. The value $t = T_N^*$, at which $\psi_N^*(T_N^*) = \Theta_N^{*'}(T_N^*) = \infty$, is called a «blow up» point.

If in (65) $T = T_N^*$, then the estimate (65) is the best possible estimate of continuous solutions to inequality (62), since the majorant in (65) is the unique continuous solution to the equation generated by (62) through replacement of sign $\leq$ by $=$.

At $N = 1$ obtain from (62) the linear integral inequality

$$|x_1^*(t)| \leq F(t) + L_1(t) \int_0^t |x_1^*(s)| ds, \quad t \in [0, T], \tag{67}$$

and the best possible estimate (65) has the form

$$|x_1^*(t)| \leq \psi_1^*(t), t \in [0, T], 0 < T < \infty, \tag{68}$$

where

$$\psi_1^*(t) = F(t) + L_1(t) e^{\int_0^t L_1(s) ds} \int_0^t F(s) e^{-\int_0^s L_1(s_1) ds_1} ds. \tag{69}$$

If, in particular, $F(t) = const > 0$, $L_1(t) = const \geq 0$, then

$$\psi_1^*(t) = F e^{L_1 t} \tag{70}$$

and (68) turns into the standard Gronwall-Bellman inequality. As is established in (Apartsyn & Spiryaev, 2010) at some values of constant parameters $F, L_m, M_m$ in (62) $\psi_N^*(t)$, $N \geq 2$, is expressed in terms of the second real branch of the Lambert function $W_{-1}(\cdot)$; in (Apartsyn, 2011) this fact is used to substantiate the stability of continuous solution to equation (2) and in (Apart-syn & Spiryaev, 2011) this result is obtained for the general case. Thus, the Lambert function or its analogs play the role of exponent in the theory of linear differential and Volterra integral equations in the case of polynomial Volterra equations.

In (Apartsyn, 2004b; Apartsyn, 2007). the authors use (34) to show that besides the unique continuous solution $x_2^*(t)$, there is the second solution $x_2^{**}(t)$, which belongs to the space of generalized functions:

$$x_2^{**}(t) = -x_1^*(t) - \frac{1}{\lambda}\delta(t),$$

where $\delta(t)$ is the Dirac $\delta$- function. The research on the structure of solutions to equation (3) in the space of generalized functions and also some partial classes of (2) is done in (Sidorov & Sidorov, 2011a; Sidorov & Sidorov, 2011b; Solodusha, 2011a; Solodusha, 2012b).

The problem of numerically solving a bilinear (quadratic) Volterra equation of the first kind was studied in (Apartsyn, 2007). Convergence of the quadrature method of right-hand rectangles on the interval $[0, \widehat{T}]$, $\widehat{T} < T_2^*$ is proved. The proving procedure is based on the Lambert function properties. The effect of the emergence of a boundary layer of the numerical solution errors in the vicinity of the "blow up" point is revealed.

In (Solodusha, 2008; Solodusha, 2009; Solodusha, 2012b) the authors considered partial classes of systems of the bilinear Volterra equations of the first kind. On the whole, however, the development of the theory and numerical methods of solving the polynomial Volterra equations are far from being completed.

The research (Solodusha, 2011b) was devoted to numerically solving similar equations related to the problem of automatic control of heat exchange processes in the case of no feedback. Some results (Solodusha, 2011c) and (Solodusha, 2012c) are presented in the final section of this paper.

## POLINOMIAL VOLTERRA EQUATIONS IN THE PROBLEM OF AUTOMATIC CONTROL

Let the input signal $x(t)$ be a vector function of time that consists of $p \geq 2$ components $x_i(t), i = \overline{1,p}$, response $y(t)$ – a scalar (which does not decrease generality) function of time, which continuously depends on $x(t)$. Then, the mathematical model of the input-output dynamic system can be represented by the Volterra polynomial of the $N$-th power:

$$y(t) = \sum_{n=1}^{N} \sum_{1 \leq i_1 \leq .. \leq i_n \leq p} f_{i_1,...,i_n}(t), \tag{71}$$

$$f_{i_1,\ldots,i_n}(t) = \int_0^t \ldots \int_0^t K_{i_1,\ldots,i_n}(t,s_1,\ldots,s_n) \prod_{m=1}^n x_{i_m}(s_m)ds_m, t \in [0,T], \tag{72}$$

where $y(t)$ – scalar function of time, $y(0) = 0, y(t) \in C^{(1)}_{[0,T]}$. The Volterra kernels $K_{i_1,\ldots,i_n}$ in (72) are symmetric only with respect to the variables that correspond to coinciding indices. Assume further that the identification problem of Volterra kernels $K_{i_1,\ldots,i_n}$ in (72) is solved. We will also suppose that the disturbances $x_i(t), i = \overline{2,p}$, in (72) are known. Along with required smoothness of the input data in (71) and (72) we will assume that $K_1(t,t) \neq 0 \ \forall t \in [0,T]$. Consider the case $N = 2$ in (71) which is the most interesting for applications. In the case of stationary dynamic system instead of (71), (72) have

$$\sum_{i=1}^p V_{1,i}x_i + \sum_{i=1}^p V_{2,i}x_i^2 + \sum_{i=2}^p \sum_{j=1}^{i-1} V_{2,ji}(x_j,x_i) = y(t), t \in [0,T], \tag{73}$$

where

$$V_{1,i}x_i \equiv \int_0^t K_i(t-s)x_i(s)ds,$$

$$V_{2,i}x_i^2 \equiv \int_0^t \int_0^t K_{ii}(t-s_1,t-s_2)x_i(s_1)x_i(s_2)ds_1 ds_2,$$

$$V_{2,ji}(x_i,x_j) \equiv \int_0^t \int_0^t K_{ji}(t-s_1,t-s_2)x_j(s_1)x_i(s_2)ds_1 ds_2, i \neq j; i,j = \overline{1,p}.$$

Consider the problem of stabilization (regulation) which is related to the search for control action $x_1(t)$ that supports the input signal $y(t)$ at a specified level $y^*$. Such a statement is possible when it concerns the problems of automatic control of technical objects. In this case equation (73) is a polynomial Volterra equation of the first kind, and its continuous solution is of a local character. In (Solodusha, 2011b) consideration is given to a numerical scheme of solving the polynomial equation (73) for $p = 2$ provided there is no feedback. Extending the research started in (Solodusha, 2011c) we will consider an algorithm for obtaining the control action $x_1(t)$ taking into account a posteriori data about deviation of the output variable $y(t)$ from the sought value $y^*$, so that $x_1(t) = u(t-h), u(\xi) = 0, \ \xi \in [-h,0], h$ – is a known constant delay. In this case the problem of the nonlinear dynamic object regulation is reduced to the search for continuous solution $u^*(t)$ to the polynomial Volterra equation of the first kind

$$V_{1,1}u + \sum_{i=2}^p V_{1,i}x_i + \sum_{i=2}^p V_{2,1i}(u,x_i) + V_{2,1}u^2 + \sum_{i=2}^p V_{2,i}x_i^2 + \sum_{i=2}^p \sum_{j=2}^{i-1} V_{2,ji}(x_j,x_i) = f(t), \tag{74}$$

where $f(t) = \varepsilon(t) - \varepsilon(t-h)$, $t \in [0,T]$.

The signal $\varepsilon(t) = y^* - y(t)$, $\varepsilon(\xi) = 0$, $\xi \in [-h, 0]$, is considered to be miscoordination or an error in control.

Assume that the numerical solution to the polynomial equation (74) exists. Find it by the cubature method of middle rectangles. Introduce a net of nodes $t_i = ih$, $t_{i-\frac{1}{2}} = \left(i - \frac{1}{2}\right)h$, $i = \overline{1,n}$, $nh = T$. Approximate the integrals in (74) by sums. In order to find the approximation $u^*(t)$ at the $\left(i - \frac{1}{2}\right)$-th node we obtain the quadratic equation with respect to $u^h_{i-\frac{1}{2}}$

$$h^2 K_{11_{\frac{1}{2},\frac{1}{2}}} \left(u^h_{i-\frac{1}{2}}\right)^2 + h\left(K_{1_{\frac{1}{2}}} + 2h \sum_{j=2}^{i} K_{11_{\frac{1}{2},k-\frac{1}{2}}} u^h_{i-j+\frac{1}{2}} + h \sum_{\mu=2}^{p} \sum_{k=1}^{i} K_{1\mu_{\frac{1}{2},k-\frac{1}{2}}} x^h_{\mu_{i-k+\frac{1}{2}}}\right) u^h_{i-\frac{1}{2}} = \quad (75)$$

$$= f(ih) - z(ih) - h \sum_{j=2}^{i} \left(K_{1_{j-\frac{1}{2}}} + h \sum_{k=2}^{i} K_{11_{j-\frac{1}{2},k-\frac{1}{2}}} u^h_{i-k+\frac{1}{2}} + h \sum_{\mu=2}^{p} \sum_{k=1}^{i} K_{1\mu_{j-\frac{1}{2},k-\frac{1}{2}}} x^h_{\mu_{i-k+\frac{1}{2}}}\right) u^h_{i-j+\frac{1}{2}},$$

where

$$z(ih) = h \sum_{\mu=2}^{p} \sum_{j=1}^{i} \left(K_{\mu_{j-\frac{1}{2}}} + h \sum_{k=1}^{i} K_{\mu\mu_{j-\frac{1}{2},k-\frac{1}{2}}} x^h_{\mu_{i-k+\frac{1}{2}}} + h \sum_{\nu=2}^{\mu-1} \sum_{k=1}^{i} K_{\nu\mu_{j-\frac{1}{2},k-\frac{1}{2}}} x^h_{\nu_{i-k+\frac{1}{2}}}\right) x^h_{\mu_{i-j+\frac{1}{2}}},$$

$$y(lh) = h \sum_{j=1}^{l} \left(K_{1_{j-\frac{1}{2}}} + h \sum_{k=1}^{l} K_{11_{j-\frac{1}{2},k-\frac{1}{2}}} u^h_{l-k-\frac{1}{2}} + h \sum_{\mu=2}^{p} \sum_{k=1}^{l} K_{1\mu_{j-\frac{1}{2},k-\frac{1}{2}}} x^h_{\mu_{l-k+\frac{1}{2}}}\right) u^h_{l-j-\frac{1}{2}} + z(lh),$$

$y(h) = z(h)$, $f(h) = \varepsilon(h)$, $f(lh) = \varepsilon(lh) - \varepsilon((l-1)h)$, $\varepsilon(ih) = y^* - y(ih)$, $l = \overline{2,n}$.

The choice of the required root in (75) is determined by the condition

$$u^h_{\frac{1}{2}} \xrightarrow[h \to 0]{} u(0) = \frac{y'(0)}{K_1(0)}.$$

Assuming in (74) $p = 2$, $u(t) = \Delta D(t)$, $x_2(t) = \Delta Q(t)$, we consider the problem of choosing $\Delta D(t)$, such that provides at specified $\Delta Q(t)$ non-zero deviation $\Delta i(t)$ from its stationary value $i_0 = 434$(kJ/kg). The scheme presented in this section was added in the software that makes it possible to identify and test the model of (71) in the case $N = 2$ as applied to the mathematical model of heat exchanger (32).

The results of computational experiments on solving the problem of automatic control in terms of feedback are presented in the Figures below. The numerical calculations of $\Delta D_1(t)$, $\Delta D_2(t)$ were made on the basis of integral models. The algorithms of their construction are described in (Apartsyn, Solodusha, Tairov & Khudyakov, 1994) and (Solodusha, Spiryaev & Scherbinin, 2007). The test signal amplitudes applied for construction of the integral models as before accounted for 25 percent of the stationary values $D_0$ and $Q_0$.

Figures 3(a) and 4(a) present the disturbances $\Delta Q(t)$ set in advance. Figures 3(b) and 4 (b) show the calculated values of solution $\Delta D_1(t)$, $\Delta D_2(t)$, which generate corresponding output signals $\Delta i_1(t)$, $\Delta i_2(t)$ at $t \in [0,30]$. It is seen that at the end of the transient process ($T = 30$ sec) $\Delta i_1(t)$ and $\Delta i_2(t)$ are equal to zero.

The main specific features of the equation of type (74) as in the scalar case relate to the local character of the region of existence of its real continuous solution. The computational experiments have shown that this property of (74) manifests itself at a disturbance $\Delta Q > 55\% Q_0$, which means potential loss of controllability of the studied heat exchange process.

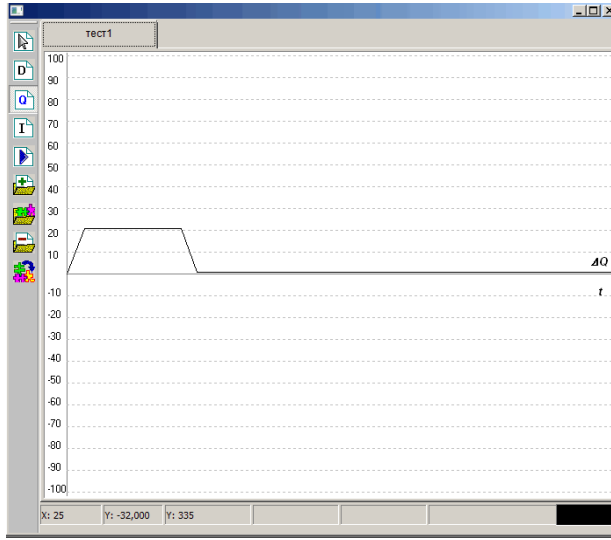

*Figure. 3(a). Disturbance $\Delta Q(t)$ (in %).*

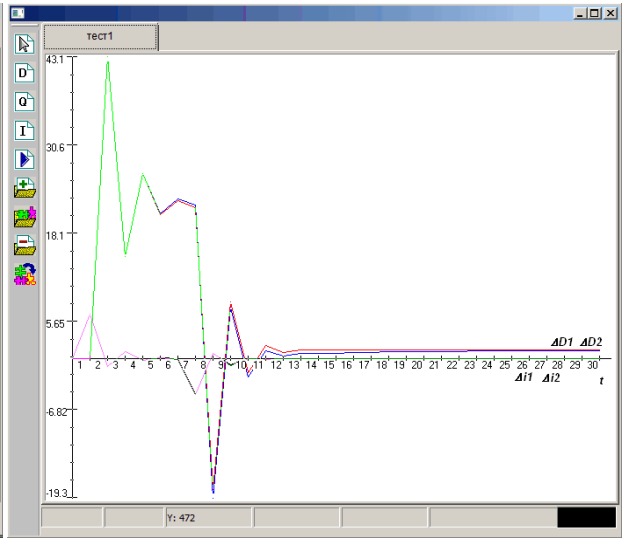

*Figure. 3(b). Control signals $\Delta D_1(t)$, $\Delta D_2(t)$ and output signals $\Delta i_1(t)$ and $\Delta i_2(t)$ (in %).*

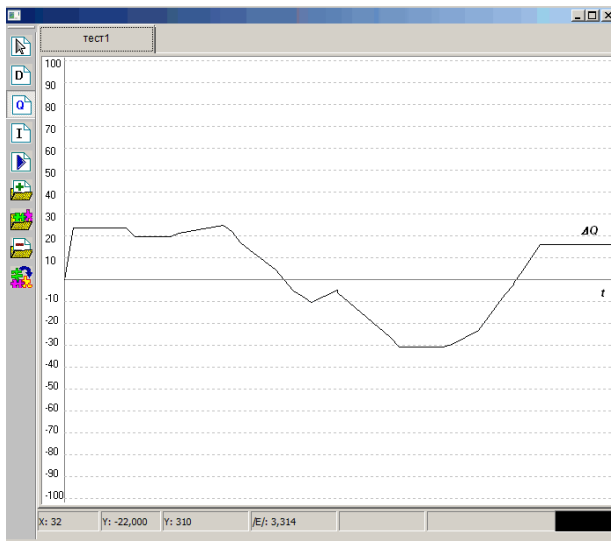

*Figure. 4(a). Disturbance $\Delta Q(t)$ (in %).*

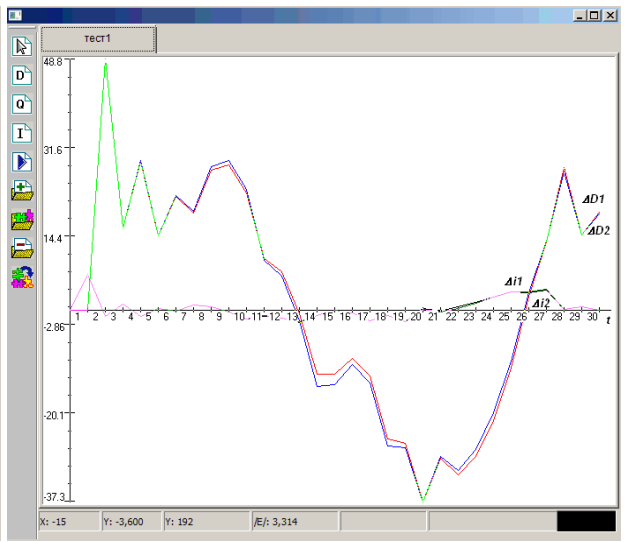

*Figure. 4 (b). Control signals $D_1(t)$, $\Delta D_2(t)$ and output signals $\Delta i_1(t)$ and $\Delta i_2(t)$ (in %).*

## CONCLUSIONS

Thus, the paper presents the main theoretical results obtained at the Energy Systems Institute SB RAS in the field of mathematical modeling of nonlinear dynamic systems with Volterra polynomials as well as examples of their application in modeling of heat exchange dynamics.

The research was supported by the grant of the RFBR No. 12–01–00722.